\documentclass[11pt,reqno]{amsart}
\usepackage{amssymb}
\newtheorem{thm}{Theorem}

\begin{document}

\bibliographystyle{plain}

\title[Rado\v s Baki\'c]{On location of zeroes of the first derivative of a polynomial}

\author{Rado\v s Baki\'c \\ \tiny{ Teacher Training Faculty, University of Belgrade}}

\address{Teacher Training Faculty, University of Belgrade, Kraljice Natalije 43, 11000 Belgrade, Serbia}
\email{\rm bakicr@gmail.com}

\begin{abstract}
Let $p(z) = z^n + a_{n-1}z^{n-1} + \cdots + a_1 z + a_0$ be a complex polynomial, $\deg p = n \geq 2$, and let $z_1, \ldots, z_n$ be the
zeros of $p(z)$, counting multiplicities. Let $c$ be the centroid of the first $n-1$ zeros $z_1, \ldots, z_{n-1}$ and let $C$ be a disc
centered at $c$ which contains $z_1, \ldots, z_{n-1}$. Then the same disc contains at least
$$\left[ \frac{n-1}{2} \right]$$
zeros of the polynomial $p'(z)$.
\end{abstract}

\maketitle

\footnotetext[1]{\, Mathematics Subject Classification 2010 Primary 26C10.  Key words
and phrases: Zeros of a polynomial, zeros of derivative, location of zeros}

\section{Introduction and statement of the result}

One of the main problems in the theory of geometry of polynomials is the following one: given some information on the location of the zeros
of a polynomial $p(z)$, deduce something on the location of one or all zeros of the polynomial $p(z)$. There are many results of this type,
we mention Gauss-Lucas theorem and Grace-Heawood theorem as famous theorems in this area. In this note we prove the following result, which deals with the same problem.

\begin{thm}
Let $p(z) = z^n + a_{n-1}z^{n-1} + \cdots + a_1 z + a_0$ be a complex polynomial, $\deg p = n \geq 2$, and let $z_1, \ldots, z_n$ be the
zeros of $p(z)$, counting multiplicities. Let
$$c = \frac{z_1 + \cdots + z_{n-1}}{n-1}$$
and let $C$ be a disc
centered at $c$ which contains $z_1, \ldots, z_{n-1}$. Then the same disc contains at least
$$\left[ \frac{n-1}{2} \right]$$
zeros of the polynomial $p'(z)$.
\end{thm}

Note that the interesting cases are $n \geq 3$, for $n=2$ there is nothing to prove.

The proof we give uses Walsh Coincidence Theorem, see \cite{Wa}, for readers' convenience we state it below.

\begin{thm}[Walsh]
Let $f(x_1, \ldots, x_n)$ be a complex polynomial, which is linear and symmetric in each of the variables $x_1, \ldots, x_n$, and has total
degree $n$. If $n$ points $x_1^0, \ldots, x_n^0$ belong to a disc $C$ (which can be either open or closed), then there is a point $x$ in $C$ such that $f(x_1^0, \ldots, x_n^0) = f(x, x, \ldots, x)$.
\end{thm}

\section{Proof}

If $z_n$ belongs to $C$, then the result follows from the Gauss-Lucas theorem. Hence we assume $z_n$ is not in $C$. Clearly, we can assume
$z_n = 0$. Moreover, using suitable rotation, we can assume $c > 0$.

Now, let $y_1, \ldots, y_{n-1}$ be the zeros of $p'(z)$, counting multiplicities. Then we have, for every $1 \leq i \leq n-1$,
$$\sum_{k=1}^n \prod_{j \not= k}(y_i - z_j) = 0.$$
The expression above is linear and symmetric in $z_1, \ldots, z_{n-1}$ and has total degree $n-1$, hence we can apply Theorem 2. Therefore
there is, for each $1 \leq i \leq n-1$, a point $c_i$ in $C$ such that
$$(n-1)y_i(y_i - c_i)^{n-2} + (y_i - c_i)^{n-1} = 0.$$
Note that we used here $z_n = 0$. From the above equation we conclude that, for each $1 \leq i \leq n-1$, we have $y_i = c_i$ or $y_i = c_i/n$. Assume the first possibility occurs $k$ times, by renumbering we can assume $y_i = c_i$ for $1 \leq i \leq k$. Clearly it suffices to show that $k >(n-3)/2$. By the above we have
$$y_i = \frac{c_i}{n}, \qquad k+1 \leq i \leq n-1.$$
Using Viet's rule, we have
$$\frac{1}{n} \sum_{i=1}^n z_i = \frac{1}{n-1} \sum_{i=1}^{n-1} y_i,$$
which gives
\begin{equation}\label{ast}
\frac{n-1}{n} c = \frac{1}{n-1} \left( \sum_{i=1}^k y_i + \sum_{i = k+1}^{n-1} \frac{c_i}{n} \right).
\end{equation}
Let $r$ be the radius of $C$. Then we have
$$ \Re \left( \frac{c_i}{n} \right) \leq \frac{r+c}{n}, \quad \Re y_j \leq r + c, \qquad 1\leq j \leq k, \quad k+1 \leq i \leq n-1.$$
Since $c > 0$, we can take real parts in (\ref{ast}) and obtain
$$\frac{n-1}{n} c \leq \frac{1}{n-1} \left( \sum_{j=1}^k (r+c) + \sum_{i=k+1}^{n-1} \frac{r+c}{n} \right) = \frac{r+c}{n-1}
\left( k + \frac{n-k-1}{n} \right),$$
which is equivalent to
$$(n-k-2) c \leq r(1+k).$$
Since $C$ does not contain the origin, we have $r < c$ and this gives the needed estimate $n-k-2 < 1 +k$.


\begin{thebibliography}{99}

\bibitem{Pr} V. Prasolov, {\it Polynomials}, Springer 2004.

\bibitem{Wa} Walsh, J. L.  {\it On the location of the roots of certain types of polynomials}, Trans. Amer. Math. Soc. 24 (1922) 163-180.




\end{thebibliography}
\end{document}